\documentclass[11pt]{article}

\usepackage{amscd}
\usepackage{amssymb,amsmath}

\newcommand{\U}{\widehat{U}}
\newcommand{\X}{\widehat{X}}
\newcommand{\J}{\mathcal{J}}
\newcommand{\m}{\mathfrak{m}}
\newcommand{\plimit}{\underleftarrow{\lim}\;}

\newtheorem{theorem}{Theorem}
\newtheorem{prop}[theorem]{Proposition}
\newtheorem{lemma}[theorem]{Lemma}
\newtheorem{cor}[theorem]{Corollary}

\date{March 6, 2008}
\author{Vladimir Baranovsky}
\title{Algebraization of bundles on non-proper schemes.}

\begin{document}

\maketitle

\section{Introduction}

This work is a contribution toward an algebraic understanding of the 
Uhlenbeck compactification. Recall, cf. \cite{DK} that for a complex
projective surface $S$ the moduli space $M_n$ of semistable vector 
bundles with fixed rank, determinant and $c_2 = n$ is non-compact, but the
union $Uhl_n = \coprod_{s \geq 0} M_{n-s} \times Sym^s S$ can be given 
a topology of a compact space (since one deals with semistable
bundles for $s \gg 0$ the space $M_{n-s}$ will be empty). We will call
$Uhl_n$ the Uhlenbeck moduli space although sometimes this name is
reserved for the closure of $M_n$ in $Uhl_n$. 

Some time ago, see e.g. \cite{L3}, \cite{BFG}, \cite{FGK}, the Uhlenbeck 
moduli space started to appear in algebraic geometry and higher 
dimensional Langlands Program. For instance, it is a convenient tool for the
study of higher versions of Hecke correspondences which modify a vector
bundle on $S$ (more generally, a principal bundle) along a divisor,
obtaining a new bundle. 
For several reasons, we would like to have a definition of $Uhl_n$ as a
``functor", i.e. we want to be able to describe in geometric terms 
the set of maps $F(T)$ (actually, a category of maps) from any test scheme
$T= Spec(A)$ to $Uhl_n$. Firstly, that would allow to define $Uhl_n$ over
any field $k$ and not to require stability. Secondly, in the study 
of the cohomology of $Uhl_n$ and the
action of Hecke correspondences on it, one needs to deal with the
phenomenon of unexpected dimension of $Uhl_n$. A possible approach 
involves defining a ``derived moduli space" $D Uhl_n$ in the sense of 
\cite{L2}  which would amount
to considering more general ``spaces" $T$. Thus, defining $Uhl_n$ as
a functor is a necessary preliminary step to constructing $DUhl_n$.

Very roughly, it is expected that a map $T\to Uhl_n$ should be described by 
a vector bundle $F$ on an open subset $U\subset T \times S$ such that
its complement $Z$ is finite over $T$, a family $\xi$ of effective 
zero cycles on $X$ parametrized by $T$ plus an agreement condition 
between $\xi$ and $F$. Such a definition gives a ``reasonable space" 
$Uhl_n$ if it satisfies a criterion due to Artin, cf. \cite{Ar}, or its 
``derived" generalization proved in \cite{L2}. 
The most difficult part of Artin's criterion is the effectiveness condition: 
if $A$ is a complete noetherian local $k$-algebra with maximal ideal 
$\mathfrak{m}$ and $A_p = A/\mathfrak{m}^{p+1}$ one needs to show that
$F(Spec(A)) = \plimit F(Spec(A_p))$. 
Ignoring the family of zero cycles $\xi$ (as will be done in this paper), 
if $X = Spec(A)\times_k S$ and $\widehat{X}$ is its formal completion
along the fiber over the closed point of $Spec(A)$, we are
trying to find whether a bundle $\mathcal{F}$ on an open subset $\widehat{U} 
\subset \widehat{X}$ comes  from a bundle $F$ on an open subset $U\subset X$. 
Such $F$ is called an \textit{algebraization} of $\mathcal{F}$.

In this paper we prove that, when $S$ has arbitrary dimension and $\widehat{U}$
has complement of codimension $\geq 3$, algebraization always exists (for
vector bundles and also for principal bundles over reductive groups). If
$\widehat{U}$ has complement of codimension $\geq 2$ then algebraization 
exists only under an additional condition (which, in the Uhlenbeck functor case, 
a guaranteed due to the presence of the relatize zero cycle $\xi$). 

Earlier similar questions were studied for coherent sheaves 
on proper schemes by Grothendieck, see \cite{EGAIII}, 
and in the case of Lefschetz type theorems by Grothendieck
and Raynaud in \cite{SGA2}, and \cite{R}. Although these results do
not apply in our case directly, our proof is based on the tools developed in 
\cite{EGAIII}, \cite{SGA2}.

In Section 2 we fix the notation, give examples illustrating some
problems to be encountered, and prove algebraization results
for vector bundles, summarized in Corollary \ref{conditions}. In Section 3
we formulate an algebraization criterion for principal bundles over
reductive groups, see Theorem 
\ref{principal}. Finally, Section 4 provides a categorical restatement of 
our results, see Theorem \ref{categorical}.

\bigskip
\noindent
\textbf{Acknowledgements.} The author thanks V. Ginzburg who first formulated
the problem of defining the Uhlenbeck functor and whose unpublished 
preprint on it (written jointly with the present author) served as a principal 
motivation for this work. Many thanks are also
due to V. Drinfeld who conjectured the statement of Theorem \ref{vector}$(i)$,  
brought the author's attention to the references \cite{SGA2} and \cite{Ha}, and 
also suggested Example 3 in Section 2.2 below. 

\medskip
\noindent
This work was supported by the Sloan Research Fellowship.

\section{Algebraization for vector bundles.}

\subsection{Setup}

We refer the reader to Expose III in \cite{SGA2} regarding basic properties
of depth and its relation to local cohomology. Let $S$ be an irreducible 
noetherian scheme of finite type over a field $k$. We will assume that $S$ 
is proper and  satisfies Serre's $S_2$ condition: for any $s \in S$, 
$depth_s \mathcal{O}_S \geq \min(\dim \mathcal{O}_{S, s}, 2)$.  
Let $V \subset S$ be an open subset with closed
complement of codimension $\geq 2$ in $S$  and
$A$ a complete noetherian local $k$-algebra with residue field $K= A/\m$ and 
associated graded $K$-algebra $gr(A) = \oplus_{p \geq 0} gr_p(A) = 
\oplus_{p \geq 0} \mathfrak{m}^p/\mathfrak{m}^{p+1}$.  Define $X = S \times_k Spec(A) $ and 
$$
X_p = S \times_k Spec(A/\m^{p+1}); \quad U_p = V \times_k Spec(A/\m^{p+1}). 
\qquad p \geq 0
$$
Let 
$
i_p: U_p \to X_p
$
be the natural open embeddings. Denote  by $\X$ the completion
of $X$ along $X_0$, which may be viewed at the limit of $\{X_p\}_{p \geq 0}$,
 cf. Section 10.6 in \cite{EGAI}. The open subset $U_0 \subset X_0$ defines an 
open formal subscheme  $\widehat{i}: \widehat{U}\to 
\widehat{X}$, given by the limit of $\{U_p\}_{p \geq 0}$. 
The ideal sheaf of $X_0$ in $X$ will be denoted by by $\J_X$ and the closed
subset $X_0 \setminus U_0$ by $Z_0$. Finally, $f:X\to Spec(A)$ is
 the natural proper projection and, for any $s \in Spec(A)$, $X_s$ 
stands for the fiber $f^{-1} (s)$.

\bigskip\noindent
Observe that $X$ may no longer satisfy the $S_2$ conditon 
(since we made no depth assumptions
on $A$). However, for $f(x) = s$ we can lift a regular sequence from 
$\mathcal{O}_{X_s, x}$ to $\mathcal{O}_{X, x}$ which gives
\begin{lemma}
\label{fiber-depth}
For any $x \in X$  with $f(x) = s$,
$
depth\; \mathcal{O}_{X, x} \geq \min (\dim \mathcal{O}_{X_s, x}, 2) . 
\qquad \square
$
\end{lemma}

\bigskip\noindent 
Consider a vector bundle $\mathcal{F}$ on $\U$, i.e. 
a sequence of vector bundles $F_p$ on $U_p$ with isomorphisms
\begin{equation}
\label{compat}
F_p|_{U_{p-1}}\simeq F_{p-1}; \qquad p \geq 1.
\end{equation}
\textbf{Definition.} We will say that a vector bundle $\mathcal{F}$ on 
$\widehat{U}$ \textit{admits an algebraization} if there exists
an open subset $U \subset X$ with $U\cap X_0 = U_0$  and a vector bundle
$F$ on $U$ such that  $\mathcal{F}$ is isomorphic to the completion of $F$, 
i.e. for $\mathcal{J}_U = \mathcal{J}_X |_U$ there exist
isomorphisms $F_p \simeq F/\mathcal{J}_U^{p+1} F$
compatible with \eqref{compat}.  
In Section 3 we apply similar terminology to principal bundles.

\medskip
\noindent
Let $Z$ be the closed subset $X \setminus U$ and 
$i: U\hookrightarrow X$ the open embedding.

\begin{lemma}
\label{fiber-codim}
Let $codim_{X_0} Z_0 \geq 2$ and suppose that $U\subset X$ is 
an open subset such that $U \cap X_0 = U_0$. For any $s \in Spec(A)$ 
define $Z_s = Z \cap X_s$. Then 
$codim_{X_s} Z_s \geq 2$ for all $s \in Spec A$ and $codim_X Z\geq 2$.
\end{lemma}
\textit{Proof.} Since $f$ is proper, the image $f(\overline{Z}_s)$
contains the unique closed point $s_0 \in Spec(A)$. Therefore  $\overline{Z}_s \cap X_0 \subset Z_0$ is not empty. By semicontinuity of 
dimensions in the fibers we have $codim_{X_s} Z_s \geq codim_{X_0}
(\overline{Z}_s \cap X_0) \geq codim_{X_0} Z_0 = 2$. 
The second assertion of the lemma follows from the first.$\square$

\bigskip
\noindent
In our discussion, we repeatedly use the following results
\begin{prop} 
\label{facts}
In the notation introduced above
\begin{enumerate}
\item Completion along $X_0$ induces an equivalence between the category 
of coherent sheaves on $X$ and the category of coherent sheaves on 
the formal scheme $\widehat{X}$.
\item For any locally free sheaf $F$ (resp. $F_0$) on $U$ (resp. $U_0$)
its direct image $i_*F$ (resp $(i_0)_* F_0$) is coherent. If
$codim_{X_0} Z_0 \geq 3$ then $R^1 (i_0)_* F_0$ is also coherent. 
\item Let $E$ be a coherent sheaf on $X$ and $\psi: E \to i_*i^* E$
the canonical morphism. Then $\psi$ is an isomorphism if and only if 
$depth_x E \geq 2$ for any point $x \in Z = X \setminus U$.
\end{enumerate}
\end{prop}
\textit{Proof.} Part (i) follows from Corollary 5.1.6 in \cite{EGAIII}.
To check the coherence of $i_* F$, by Corollary VIII.2.3 in \cite{SGA2}
it suffices to check that $depth_x F \geq 1$ for any point $x \in U$ such 
that $\overline{\{x\}} \cap Z$ has codimension 1 in $\overline{\{x\}}$.
But Lemma 1 and local freeness of $F$  imply that any $x$ with $depth_x F = 0$
must be generic in its fiber, and the Lemma 2 implies that $\overline{\{x\}} 
\cap Z$ would in fact have codimension 2 in $\overline{\{x\}}$. The same
proof applies to $(i_0)_* F_0$. If $codim_{X_0} Z_0 \geq 3$ then the
above argument can also by applied to $R^1 (i_0)_* F_0$ once we show that
$depth_x F_0 \geq 2$  for any $x \in U_0$ such that  
$\overline{\{x\}} \cap Z_0$ has codimension 1 in $\overline{\{x\}}$. But
by $S_2$ condition $depth_x F_0 \leq 1$ can only hold for points $x$ of
codimension $\leq 1$ in $U_0$, which would imply that $\overline{\{x\}} 
\cap Z_0$ has codimension $\geq 2$ in $\overline{\{x\}}$. This proves $(ii)$. 
Part $(iii)$ is a particular case of Corollary II.3.5 in \textit{loc.cit}.
$\square$

\subsection{Examples.}

The first example with $codim_{X_0} Z_0 = 3$ and $K = k$
shows that one may not be able to take $U= U_0 \times_k Spec(A)$. 

\bigskip
\noindent
\textbf{Example 1.}
Take $S = X_0 = \mathbb{P}^3$ with homogeneous coordinates $[x:y:z:w]$
and set $V = U_0 = S \setminus [0:0:0:1]$,  $A= k[[t]]$ 
(formal power series in $t$).
Define vector bundles $F_p$ as kernels of 
$$
\varphi_p: \mathcal{O}_{U_p}^{\oplus 3} \to \mathcal{O}(1)_{_{U_p}}; \quad 
(s_1 \oplus s_2 \oplus s_3) \mapsto s_1 x + s_2 y + s_3 (z-tw).
$$
Observe that $\varphi_p$ is surjective since $t$ is nilpotent on $U_p$ and
$[0:0:0:1] \notin U_p$. 
\begin{lemma} The bundle $\mathcal{F}$ admits no algebraization $(U, F)$ 
with  $U= U_0
\times_k Spec(A)$. 
\end{lemma}
\textit{Proof.}
Set $F$ to be the kernel of morphism $\varphi: \mathcal{O}^{\oplus 3} \to 
\mathcal{O}(1) $of vector bundles on $U$, given by the same formula as for 
$\varphi_p$.  By definition, $\varphi$ is not surjective only at
 $P= [0: 0: t: 1] \in U$ which projects to the
 generic point $\xi = Spec(k[t^{-1},t]]) \in Spec(A)$. The specialization 
at $t = 0$ is not  in $U_0$, hence $P$ is closed in 
$U$ and $U\setminus P$ is an open subset containing $U_0$. Since on  
$U\setminus P$ we have the short exact 
sequence of locally free sheaves
$$
0 \to F \to \mathcal{O}^{\oplus 3} \to \mathcal{O}(1) \to 0,
$$
the restriction of $F$ to each  $U_p$ is given by  $F_p$, i.e. $\mathcal{F}$ is
indeed the completion of $F$. 
On the other hand, $F$ is not locally free at $P$: from 
$0 \to F \to \mathcal{O}^{\oplus 3}_U \to \mathcal{O}_U \to k_P \to 0$
we immediately get $\mathcal{E}xt^1(F, \mathcal{O}_U) \simeq \mathcal{E}xt^3(k_P, \mathcal{O}) \simeq k_P$ since the middle two terms
are projective.

Suppose that $E$ is a locally free sheaf on $U$ with completion 
isomorphic to $\mathcal{F}$. We will see later in Proposition 
\ref{converse}$(ii)$ that in such situation we must have:
$
\widehat{i_* E} \simeq \widehat{i}_* \mathcal{F} \simeq \widehat{i_* F}
$
hence by Proposition \ref{facts}$(i)$,  $i_* F \simeq i_* E$ which contradicts 
$\mathcal{E}xt^1(F, \mathcal{O}_U) \neq 0$.  $\square$

\bigskip
\noindent
The second example illustrates that for $codim_{X_0}\; Z_0 = 2$, a pair $(U, F)$
may not exist at all.

\bigskip
\noindent
\textbf{Example 2.}
Consider $A = k[[t]]$ and
$S = X_0 = \mathbb{P}^2$ with homogeneous
coordinates $(x:y:z)$. Let $V = U_0 = X_0 \setminus P$
where $P= (0:0:1)$  and define a rank 2 bundle $F_p$
on $U_p = U_0 \times_k Spec(k[t]/t^{p+1})$ as follows. The affine open
subsets $U^{(x)}_p$, $U^{(y)}_p$ given by non-vanishing of $x$, resp.
$y$, form a covering of $U_p$ and we can glue trivial rank 2
bundles on these open sets, using the transition function 
$$
\left( \begin{array}{cc}
1 &  \sum_{m=0}^p \big(\frac{t z^2}{xy}\big)^m \\
0  & 1 \end{array} \right)
$$
on $U^{(x)}_p \cap U_p^{(y)}$. 
Clearly $F_p|_{U_{p-1}} \simeq F_{p-1}$ in a natural way, 
and we obtain a 
vector bundle $\mathcal{F}$ on $\widehat{U}$. 
\begin{lemma}
There exists no vector bundle $F$ on $U = X \setminus Z$ with 
$\widehat{F} \simeq \mathcal{F}|_{\widehat{U} \setminus (Z \cap U_0)}$, 
for any closed subset $Z \subset X$ such that $Z_0 \subset (Z \cap X_0)$ and
$codim_{X_0} (Z\cap X_0)\geq 2$.  
\end{lemma}
\textit{Proof.}
Suppose otherwise and take the direct image of $F$ with 
respect to the open embedding $i: U \to X$.  By Proposition \ref{facts}, 
$i_* F$ is coherent and has $depth \geq 2$ at all codimension $2$ points of $X$. 
Since modules of depth 2 over two-dimensional regular local 
rings are free by Auslander-Buchsbaum formula,  $i_*F$ will 
be locally free in codimension two. Therefore 
shrinking $Z$ we can assume that $Z$ has
codimension 3 in $X$ which in our case means that $Z$ is 
 a finite set of points in $X_0$.  
Then the short exact sequence of sheaves on 
$X \setminus Z$ 
$$
0 \to F \stackrel{t^{p+1}}\longrightarrow F \to F_p \to 0,
$$
leads to a long exact 
sequence on $X$:
$$
0 \to i_*F \stackrel{t^{p+1}}\longrightarrow i_*F \to (i_p)_* F_p \to
R^1i_*F \stackrel{t^{p+1}}\longrightarrow R^1 i_*F
$$
where $R^1i_* F$ is coherent for the same reason as in Proposition 
\ref{facts}$(ii)$. Since $R^1i_* F$ is supported at the finite set $Z$ of 
closed points, it has finite length at each of them and the last
arrow is zero for $p \geq p_0$. For such $p$ we can write 
$
i_*F \to (i_p)_* F_p \to R^1i_*F \to 0 
$
which gives
$$
i_*F \otimes_{\mathcal{O}_X} k(P) \to 
(i_p)_* F_p \otimes_{\mathcal{O}_X} k(P) 
\to R^1i_*F \otimes_{\mathcal{O}_X} k(P) \to 0
$$
To prove the lemma  it suffices to show 
that $\dim_{k} (i_p)_* F_p \otimes_{\mathcal{O}_X} k(P)$
is unbounded as $p \to \infty$. 

\bigskip
\noindent
To that end, replace $X_0$ with the affine open 
subset $\widetilde{X}_0 \simeq \mathbb{A}^2$ given by non-vanishing of
$z$, with affine coordinates $ u = \frac{x}{z}$, $v = \frac{y}{z}$. 
Set $W_0 = U_0 \cap \widetilde{X}_0$ and similarly for $\widetilde{X}_p$, $W_p$, 
$W_p^{(x)}$ and $W_p^{(y)}$. Then $(i_p)_* F_p|_{\widetilde{X}_p}$ is the sheaf
associated to $H^0(W_p, F_p|_{W_p})$
viewed as a module over $A(\widetilde{X}_p) =  k[u, v, t]/t^{p+1}$. 
By its definition, $F_p$ is an extension of $\mathcal{O}_{U_p}$
with $\mathcal{O}_{U_p}$
which leads to long exact sequence
$$
0 \to H^0(W_p, \mathcal{O}_{W_p})\to 
H^0(W_p, F_p|_{W_p}) \to 
H^0(W_p, \mathcal{O}_{W_p})\to 
H^1(W_p, \mathcal{O}_{W_p}).
$$
where the last arrow sends the constant function $1$ to the class of the 
extension. 
Let $M_p$ be the kernel of the last arrow. It suffices to 
show that $dim_k (M_p/\langle u, v, t \rangle M_p)$ is unbounded. 
Computing  $M_p$ via the affine covering $\{W^{(x)}_p$, 
 $W^{(y)}_p\}$ we identify it with the kernel of 
$$
k[u, v, t]/t^{p+1} \stackrel{\pi_p \circ \psi_p}
\longrightarrow \frac{1}{uv}k[u^{-1}, v^{-1}, t]/t^{p+1}
$$
where $\psi_p$ is  multiplication by 
$\sum_{l=0}^p\big(\frac{t}{uv})^l$
(i.e. the upper right corner of the transition matrix 
in the definition of $F_p$), and $\pi_p$ is the natural 
projection 
$$
k[u, u^{-1}, v, v^{-1}, t]/t^{p+1} \to 
\frac{1}{uv}k[u^{-1}, v^{-1}, t]/t^{p+1}
$$ 
It follows that  $M_p$ is 
generated by the monomials $t^p, t^i u^{p-i}, t^iv^{p-i}$
for $i = 0, \ldots, p-1$, thus
$$
\dim_k (M_p/\langle u, v, t \rangle M_p) = 2p+1 \to \infty \qquad
\textrm{as\ } p \to \infty. \qquad \square
$$ 

\bigskip
\noindent
\textbf{Example 3.} (Suggested to the author
by V. Drinfeld.) The bundle in the previous example has trivial determinant,
but if we don't insist on this condition, then there is a rank one example:
glue two trivial bundles on 
$U^{(x)}_p, U^{(y)}_p$ using the transition function 
$\sum_{m = 0}^p \big( \frac{t z^2}{xy}\big)^m$.
The resulting line bundle admits no algebraization since again
$\dim_k (i_p)_* F_p \otimes_{\mathcal{O}_X} k(P)$ is not bounded as
$p \to \infty$. 

\subsection{Algebraization of vector bundles.}

\begin{theorem}\label{vector} In the notation of section 2.1,
\begin{enumerate}
\item If $codim_{X_0} Z_0 \geq 3$ then $\mathcal{F}$ admits an
algebraization.

\item If  $codim_{X_0} Z_0 \geq 2$ and
the cokernel of the natural morphism $(i_p)_* F_p|_{X_{p-1}} \to 
(i_{p-1})_* F_{p-1}$ is supported in codimension $\geq 3$ for all $p$ large 
enough, then $\mathcal{F}$ admits an algebraization. 

\item In either of the two situations (codimension $\geq 3$ or
codimension $\geq 2$ with the additional support assumption)
the projective system $\{(i_p)_* F_p\}_{p \geq 0}$ 
satisfies the Mittag-Leffler condition, the direct image $\widehat{i}_*
\mathcal{F}$ is coherent and isomorphic to 
$\plimit (i_p)_* F_p$.   
\end{enumerate}
\end{theorem}
\textit{Proof.} We split the proof of $(i)$ and $(ii)$ in a number of steps.  
Part $(iii)$ will follow from Step 2. 

\medskip
\noindent
\textit{Step 1.}

\noindent
Suppose that $\widehat{i}_* \mathcal{F}$ is coherent. By 
Proposition \ref{facts}$(i)$  there exists a unique coherent
sheaf $E$ on $X$ such that $\widehat{E} \simeq \widehat{i}_* \mathcal{F}$. 
The subset $U \subset X$ of points where $E$ is locally free is open and 
contains $U_0$ (e.g. by Nakayama's Lemma). Shrinking $U$ if necessary we can achieve $U \cap X_0 = U_0$. Now set $F  = E|_U$.

\medskip
\noindent
\textit{Step 2.}

\noindent
Therefore $(i)$ and $(ii)$ are reduced to showing that, under the conditions
stated, 
$\widehat{i}_* \mathcal{F}$ is coherent. To that end we modify the argument of 
0.13.7.7 in \cite{EGAIII} which will also prove $(iii)$. 
First, as in 0.13.7.2 of \textit{loc. cit.}, 
we choose injective resolutions $F_k \to L^\bullet_k$ such 
that $L^\bullet_{k+1}/\J^{k+1}_U L^\bullet_{k+1} \simeq L^\bullet_k$
and the natural filtrations by $\J_U^n(\ldots)$ agree
with those on $F_k$. 
Each $\widehat{i}_* (L_k^\bullet)$ is a filtered complex and has a
spectral sequence with $E_1$ term given by 
$$
E_1^{pq} = R^{p+q}\widehat{i}_* (\mathcal{J}_U^p F_k/ \mathcal{J}^{p+1}_U F_k)
$$
As in 0.13.7.3 of \textit{loc.cit.}  we pass to the limit as 
$k \to  \infty$ and get a spectral sequence with 
$$
E^{p, q}_1 = R^{p+q}\widehat{i}_* (F_p/F_{p+1}) 
\simeq R^{p+q}\widehat{i}_* (F_0) \otimes_K (\mathfrak{m}^p/
\mathfrak{m}^{p+1}) =R^{p+q}\widehat{i}_* (F_0) \otimes_K gr_p(A)
$$
We are interested in the components
$$
E^0_1 = \bigoplus_{p+q =0} E^{p, q}_1 = 
\widehat{i}_* (F_0) \otimes_K gr(A); \quad
E^1_1 = \bigoplus_{p+q =1} E^{p, q}_1 = 
R^1\widehat{i}_* (F_0) \otimes_K gr(A).
$$
We would like to show that the spectral sequence converges at the 
$E^0 = \oplus E^{p, -p}$ terms. Note that each 
$E_{k+1}^1 = \oplus E^{p, 1-p}_{k+1}$ is
a quotient of a subsheaf in $E^1_k$ while each $E^0_{k+1}$ is a subsheaf 
$E^0_k$ (since $E^{p, -1 -p}$ terms are zero). Taking successive preimages 
of the boundaries in $E_{r-1}$, $E_{r-2}, \ldots, E_1$ we get a sequence 
of boundary
subsheaves $B_1 \subset B_2 \subset B_3 \subset \ldots \subset E_1^1$, 
and taking preimages of cycles in  $E_k$ we get a sequence of cycle subsheaves 
$ E^0_1 \supset Z_1 \supset Z_2 \supset Z_3 \supset \ldots .$ By 0.13.7.6
in \textit{loc.cit.}  these are actually 
$\mathcal{O}_{X_0} \otimes_K gr(A)$-submodules.

Suppose that sequence of cycles stabilizes, i.e. for some $r_0$ one has
$Z_r = Z_{r_0}$ whenever $r \geq r_0$,  
then by  0.13.7.4 in \cite{EGAIII},
 the projective system $\{\widehat{i}_* (F_k)\}_{k \geq 0}$ satisfies the 
Mittag-Leffler condition and the associated graded of
$\widehat{i}_*(\mathcal{F})$ is precisely 
$Z_{r_0} \subset  \widehat{i}_* (F_0) \otimes_K gr(A)$. 
But $\widehat{i}_*(F_0)$ is a coherent by Proposition \ref{facts}$(ii)$, 
hence the subsheaf
$gr(\widehat{i}_* \mathcal{F}) \subset  \widehat{i}_* (F_0) \otimes_K gr(A)$
is a coherent $\mathcal{O}_{X_0} \otimes_K gr(A)$-module, by the 
noetherian property of $X_0$ and $A$. By \textit{loc.cit.} 13.7.7.2,
$\widehat{i}_*\mathcal{F}$ is itself coherent on $\widehat{X}$. Also, 
$\widehat{i}_* \mathcal{F} \simeq \plimit (i_p)_* F_p$ by 0.13.7.5.1 in 
\textit{loc.cit.}.

\medskip
\noindent
\textit{Step 3.} 

\noindent
Now the assertion of the theorem is reduced to showing that the sequence
of cycles $Z_1 \supset Z_2 \supset \ldots$ stabilizes. By definition
of $Z_i$ this is equivalent to saying that the higher
differentials of the spectral seqence $d_r: E^0_r \to E^1_r$ become zero
for $r \geq r_0$. That in turn is equivalent to saying that
the sequence of boundaries $B_1 \subset B_2 \subset B_3 \subset \ldots $, 
also stabilizes. 

If $codim_{X_0} Z_0 \geq 3$ by Proposition \ref{facts}$(ii)$, 
$R^1(i_0)_* F_0$ is also coherent and $\{B_r\}_{r \geq 1}$ stabilizes by the
noetherian property of $R^1 (i_0)_* F_0 \otimes_K gr(A)$, which proves
 $(i)$. If $codim_{X_0} Z_0 \geq 2$  we need to find
a coherent subsheaf of $R^1 (i_0)_* F_0 \otimes_K gr(A)$ containing
$B_r$ for all $r \geq 1$.

\medskip
\noindent
\textit{Step 4.} 

\noindent
At this point we reduced $(ii)$ 
to showing that, under the assumptions stated, there exists a coherent
subsheaf 
$G \subset R^1 (i_0)_* F_0$ such that $B_r \subset G \otimes_K gr(A)$ for 
all $r$. 
By 0.11.2.2 in \cite{EGAIII} for $r \geq p$ the term $B^{p, 1 - p}_r$
is the image of the connecting homomorphism 
$$
\widehat{i}_* F_p \to \widehat{i}_* F_{p-1} 
\stackrel{\rho_p} \to R^1 \widehat{i}_* F_0 \otimes_K 
(\mathfrak{m}^p/\mathfrak{m}^{p+1})
$$
in the long exact sequence obtained by applying $R \widehat{i}_*$ to
the short exact sequence on $\widehat{U}$:
$$
0 \to F_0 \otimes_K (\mathfrak{m}^p/\mathfrak{m}^{p+1})
\to F_p \to F_{p-1} \to 0.
$$
Observe that by our assumptions each $Im(\rho_p)$ is coherent, and
supported in codimension $\geq 3$ for $p \gg 0$. Thefore we are done once we 
show that the subsheaf of $R^1 (i_0)_* F_0$ formed by all
sections with support in codimension $\geq 3$, is coherent whenever
$codim_{X_0} Z_0 \geq 2$ and $F_0$ is locally free on $U_0$. 

\medskip
\noindent
\textit{Step 5.}

\noindent
Set $Q = (i_0)_* F_0$, a coherent sheaf on $X_0$ by Step 2.
By the standard exact sequence we have 
$\mathcal{H}^2_{Z_0} Q = R^1 (i_0)_* Q|_{U_0}
= R^1 (i_0)_* F_0$, so it suffices to show that $\mathcal{H}^0_{\geq 3}
\mathcal{H}^2_{Z_0} Q$ is coherent where $\mathcal{H}^0_{\geq 3}$ is the functor of sections supported in codimension $\geq 3$. Let 
$\mathcal{H}^i_{\geq 3}$ be the higher  derived functors.

First, the standard spectral sequence for the composition of 
functors $R\mathcal{H}^0_{\geq 3}$, $R\mathcal{H}^0_{Z_0}$ has 
$E_2^{p, q} = \mathcal{H}^p_{\geq 3} \mathcal{H}^q_{Z_0} Q$. 
But $\mathcal{H}^i_{Z_0} Q = 0$ for $i = 0, 1$ 
by Proposition \ref{facts}$(iii)$,
so 
$$
\mathcal{H}^0_{\geq 3}\mathcal{H}^2_{Z_0} Q \simeq  
\mathcal{H}^2_\Phi Q
$$
where $\mathcal{H}^2_\Phi$ is the local cohomology with the family of supports $\Phi$ formed by all codimension $\geq 3$ closed subsets in $Z_0$. 

\medskip
\noindent
\textit{Step 6.}

\noindent
To show that $\mathcal{H}^2_\Phi Q$ is coherent let $\omega$ be
the dualizing complex of $X_0$, cf. \cite{RD}. More precisely, by 
\textit{loc.cit} $\omega$ is quasi-isomorphic to a complex of
injective sheaves
$$
0 \to \mathcal{K}^0 \to \ldots \to \mathcal{K}^{\dim_K X_0} \to 0
$$
with 
$
\mathcal{K}^i = \bigoplus_{codim_{X_0} x = i} i^x_* (I(x))
$
and $i^x: Spec(\mathcal{O}_{X_0, x}) \to X_0$ is the natural morphism,
while $I(x)$ is an injective envelope of the residue field $k(x)$
as a module over $\mathcal{O}_{X_0, x}$. An easy but important observation
which we use below, is that  $\mathcal{K}^p$ has no sections supported
in codimension $\geq p +1$. 

By definition of a dualizing complex, the double complex
$\mathcal{K}^{p, q} = \mathcal{H}om(\mathcal{H}om(Q, \mathcal{K}^{q}),
\mathcal{K}^p)$ has total complex quasi-isomorphic to $Q$. Moreover, by 
Proposition IV.2.1 and the remark on page 123 in \cite{RD}, this total
complex is a flasque resolution of $Q$ and hence can be used to compute 
$\mathcal{H}^\bullet_{\Phi}(Q)$. This leads to a spectral sequence:
$$
E_2^{p, q} = \mathcal{E}xt^{p}_{\Phi} (\mathcal{E}xt^{-q} (Q, \omega), 
\omega) \Rightarrow \mathcal{H}^{p+q}_{\Phi}(Q)
$$
where $\mathcal{E}xt^p_\Phi  = R^p (\Gamma_{\Phi} \circ \mathcal{H}om)$ 
and the $\mathcal{E}xt$ sheaves are understood in the 
sense of hypercohomology. 

Since $E_2^{p, q} \neq 0$ only for when $p$ and $(-q)$ are between 
$0$ and $\dim_K X_0$ only finitely many terms with fixed $p + q$ will be 
non-trivial and to show that $\mathcal{H}^2_\Phi(Q)$ is coherent it
suffices to show that $\mathcal{E}xt^{p}_{\Phi} (\mathcal{E}xt^{p-2} (Q, 
\omega), \omega)$ is coherent for $p \geq 2$.

\bigskip
\noindent
\textit{Step 7.}

\noindent
First observe that $\mathcal{E}xt^2_\Phi(G, \omega) = 0$ for any 
quasi-coherent sheaf $G$ since $\mathcal{K}^2$ has no sections supported
in codimension $\geq 3$ and hence no sections with support in $\Phi$. 
Hence we can assume that $p \geq 3$. 

Denote $\mathcal{R}^p = \mathcal{E}xt^{p-2} (Q, \omega)$. We first 
claim that $codim_{X_0} Supp(\mathcal{R}^p) = d \geq p \geq 3$. In fact, 
let $x \in Supp(\mathcal{R}^p)$ be a point with $\dim \mathcal{O}_{X_0, x} = d$. 
Since the stalk of $\mathcal{R}^p_x$ is non zero, and by local duality, 
cf. V.6 in \cite{RD} its completion is dual to $\mathcal{H}^{d + 2 - p}_x (Q)$
we conclude that $\mathcal{H}^{d + 2 - p}_x (Q) \neq 0$. Then $d+2 - p \geq 0$
and $d \geq p - 2 \geq 1$. If $d =  1$ then $p = 3$ and also 
$x \notin Z_0$ hence the stalk $Q_x$ is free.  Thus 
$\mathcal{H}^0_x(\mathcal{O}) \neq 0$, contradicting the $S_2$ assumption.
If $d \geq 2$ then applying the $S_2$ condition when $x \notin Z_0$ and
Proposition \ref{facts}$(iii)$ when $x \in Z_0$ we actually have $d + 2 - p
\geq 2$ so $d \geq p$ as required. 

By primary decomposition, the coherent sheaf $\mathcal{R}^p$ admits a finite
filtration by coherent subsheaves such that all successive quotients
have irreducible supports of codimension $\geq p$. By the standard long 
exact sequence for $\mathcal{E}xt^\bullet_\Phi(\cdot, \omega)$ is 
suffices to show that $\mathcal{E}xt^p_\Phi(G, \omega)$ is coherent
whenever $p \geq 3$ and $G$ is a coherent sheaf with irredicuble support
$Y$ of codimension $\geq p$. 

If $Y \nsubseteq Z_0$ for any $W$ in the family $\Phi$, the intersection 
$Y\cap W$ is not equal to $Y$ and therefore has codimension 
$\geq p +1$. But then $\mathcal{E}xt^p_\Phi(G, \omega) = 0$ because
any section 
$\rho$ of $\mathcal{H}om(G, \mathcal{K}^p)$ representing a class in 
$\mathcal{E}xt^p_\Phi(G, \omega)$ has zero values since $\mathcal{K}^p$
has no sections supported in codimension $\geq p +1$. 
If $Y \subseteq Z_0$ then $Y$ is an element of $\Phi$ and  
$\mathcal{E}xt^p_\Phi(G, \omega) \simeq 
\mathcal{E}xt^p (G, \omega)$ since  all sections of 
$\mathcal{H}om(G, \mathcal{K}^t)$ have support in $\Phi$. 
But $\mathcal{E}xt^p (G, \omega)$ is coherent which
finishes the proof.  
 $\square$

\bigskip
\noindent
The converse to Theorem \ref{vector} can be formulated as follows. 

\begin{prop} 
\label{converse}
In the setting of Section 2.1, assume that $\mathcal{F}$
admits an algebraization $(U, F)$ and view each $F_p$ as a sheaf on $U$.  
Then
\begin{enumerate}
\item The cokernel of $i_* F_p \to i_* F_{p-1}$ is supported in codimension 
$\geq 3$ for $p \gg 0$.  

\item
The isomorphism $\widehat{F} \simeq \mathcal{F}$ extends to 
direct images: $\widehat{i_* F} \simeq \widehat{i}_* \mathcal{F}$.
In particular, $\widehat{i}_* \mathcal{F}$ is coherent. 
\end{enumerate}
\end{prop}
\textit{Proof.} 
To prove $(i)$ observe that the cokernel of $i_* F_p 
\to i_* F_{p-1}$ is is annihilated by $\mathcal{J}_X$, being
a subsheaf of $R^1 i_* F_0 \otimes_K gr_p(A)$, and is therefore
isomorphic to the cokernel  of $i_* F_p|_{X_0} \to i_* F_{p-1}|_{X_0}$.  

We will first show that the natural map $i_* F_p|_{X_0} \to i_* F_0$ is
an embedding of sheaves for all $p$. Considering the exact sequence
$$
0 \to \mathcal{J}_X (i_* F_p) \to i_*F_p \to i_* F_p |_{X_0} \to 0 
$$ 
and its map to the first terms of the sequence
$$
0 \to i_* (\mathcal{J}_U F_p) \to 
i_*  F_p \to i_* F_0 \to R^1 i_* (\mathcal{J}_U F_p) \to \ldots 
$$
we see that $i_*F_p|_{X_0} \to i_*F_0$ is an embedding precisely 
when the natural map $\mathcal{J}_X (i_* F_p) \to i_* (\mathcal{J}_U 
F_p)$ is an isomorphism. Observe that $i_* \mathcal{O}_U = \mathcal{O}_X$
hence $i_* \mathcal{J}_U$ is a sheaf of ideals in $\mathcal{O}_X$. 

Using Lemma 1 and the Cohen-Macaulay assumption on $X_0$ we see  
that $\mathcal{H}^{t}_Z
\mathcal{O}_X = \mathcal{H}^t_{Z_0} \mathcal{O}_{X_0} = 0$ for $t = 0, 1$.
By the short exact sequence $0 \to \mathcal{J}_X \to \mathcal{O}_X \to \mathcal{O}_{X_0} \to 0$ we derive $\mathcal{H}_Z^t \mathcal{J}_X= 0$ for
$t = 0,1$ and hence $\mathcal{J}_X= i_*\mathcal{J}_U$ by Proposition \ref{facts}
(iii). Then 
$$
i_*(\mathcal{J}_UF_p) = (i_* \mathcal{J}_U) (i_* F_p) 
= \mathcal{J}_Xi_* F_p 
$$
as required. Similarly, $i_*F|_{X_0} \to i_*F_0$ is an embedding. So for any 
$p \geq 1$ we have embeddings 
$$
i_*F|_{X_0} \hookrightarrow i_* F_p |_{X_0} \hookrightarrow
 i_* F_{p-1}|_{X_0} \hookrightarrow i_*F_0
$$
Consequently, the coherent sheaf
$\mathcal{K} = Coker (i_*(F)|_{X_0} \to i_* F_0)$ has a decreasing filtration 
by images of $i_* F_p|_{X_0}$ and each $Coker( i_* F_p |_{X_0} \to i_* 
F_{p-1}|_{X_0})$ is its successive quotient. 
But $\mathcal{K}$ is a coherent sheaf with $Supp(\mathcal{K}) \subset Z_0$
and $Z_0$ has at most finitely many points of codimension 2.
Since for each point $x \in X_0$ of codimension 2, the localization 
$\mathcal{K}_x$ is a module of finite length, only finitely many successive 
quotients of the filtration of $\mathcal{K}$ can be non-trivial in codimension 
2, which proves $(i)$. 

To prove $(ii)$ first observe that $\widehat{i}_* \mathcal{F}$ and 
$E = i_*F$ are coherent by Theorem \ref{vector}$(iii)$ and 
Proposition \ref{facts}$(ii)$, respectively. 
By Proposition \ref{facts}$(i)$ we can
find a sheaf $E'$ such that $\widehat{E'} \simeq \widehat{i}_* \mathcal{F}$. 
The isomorphism $\widehat{E}|_{\widehat{U}} \simeq \mathcal{F} = 
\widehat{i}^* \widehat{i}_*\mathcal{F}$ extends uniquely to  a morphism of 
sheaves  $\widehat{\phi}:\widehat{E} \to \widehat{i}_* \mathcal{F} = 
\widehat{E'}$.
By Proposition \ref{facts}$(i)$, $\widehat{\phi}$ is the completion of 
a unique morphism $\phi: E\to E'$ which by Corollary 10.8.14 in 
\cite{EGAI} should be an isomorphism 
on an open subset $W$ containing $U_0$. Shrinking $W$ if necessary we 
can assume $W \subset U$. By Lemma 2, each point $x \in U\setminus W$ has
codimension $\geq 2$ in its fiber, hence $depth_x E \geq 2$ by Lemma 1.
For $x \in X\setminus U$ we still have $depth_x E \geq 2$ by 
Proposition \ref{facts}$(iii)$. Applying the same result to 
$j: W \hookrightarrow X$  instead of $U$ we see that $E = j_* j^* E$.  
By adjunction of $j^*$ and $j_*$ the isomorphism 
$(\phi|_W)^{-1}: j^* E' \to j^* E$ extends uniquely to 
a morphism $\psi: E' \to j_*j^* E = E$. 

By construction, the composition $\psi \phi: E\to E$ restricts to identity
on $W$ hence $\psi \phi = Id_E$, by the same adjunction. 
Similarly, the composition $\widehat{\phi} \widehat{\psi} =
\widehat{E'} \to \widehat{E'}$ restricts to identity on $\widehat{U}$
and since $\widehat{E'} \simeq  \widehat{i}_* \mathcal{F}$, 
we must have $\widehat{\phi} \widehat{\psi}
= Id_{\widehat{E'}}$, so $\phi \psi = Id_{E'}$ by Proposition 
\ref{facts}$(i)$. 
We have proved that $E= i_*F \simeq E'$. Since 
$\widehat{E'} = \widehat{i}_*\mathcal{F}$ we conclude that $\widehat{i_* F}
= \widehat{i}_* \mathcal{F}$. 
 $\square$

\begin{cor}
The following conditions are equivalent:
\begin{enumerate}
\label{conditions}
\item The cokernel of $(i_p)_* F_p \to (i_{p-1})_* F_{p-1}$ is supported
in codimension $\geq 3$ for $p \gg 0$.

\item The projective system $\{\widehat{i}_* F_p\}_{p\geq 1}$ satisfies the 
Mittag-Leffler condition.

\item The direct image $\widehat{i}_* \mathcal{F}$ is coherent.

\item The bundle $\mathcal{F}$ admits an algebraization.
\end{enumerate}
\end{cor}
\textit{Proof.} The implications $(i) \Rightarrow (ii)$ and 
$(iii) \Rightarrow (iv)$ are established in the proof of Theorem \ref{vector}. 
The implication $(iv) \Rightarrow (i)$ is
proved in Proposition \ref{converse}. If the projective system 
$\{\widehat{i}_* F_p\}_{p\geq 1}$ satisfies the Mittag-Leffler condition, by 
0.13.3.1 in \cite{EGAIII} the natural 
map $\widehat{i}_*\mathcal{F} \to \plimit \widehat{i}_* F_p$ is an isomorphism.
By the Mittag-Leffler condition we can replace $\widehat{i}_*F_p$ by a
system of subsheaves $G_p \subset \widehat{i}_* F_p$ so that the property 
$\widehat{i}_* \mathcal{F} \simeq \plimit G_p$ still holds and $G_p|_{X_{p-1}}
\to G_{p-1}$ is surjective. Since each $G_p$ is coherent by the noetherian
property of $X_p$, Proposition 10.11.3 in \cite{EGAI} tells that 
$\plimit G_p$ is also coherent. Therefore,  $(ii) 
\Rightarrow (iii)$. $\square$ 

\bigskip
\noindent
\textbf{Remark.} Suppose that $X_0$ is a smooth projective surface over $K$,
$\xi = k_1 P_1 + \ldots + k_l P_l$ an effective zero cycle and
$F_0$ a rank $n$  vector bundle on $U_0 = X_0 \setminus \{P_1, \ldots, P_l\}$.
The pair $(F_0, \xi_0)$ should define a point $Spec(K) \to Uhl_n$
of the Uhlenbeck functor. Assume that $(F, \xi): Spec(A)\to Uhl_n $ extends 
$(F_0, \xi_0)$. Then it is expected that $Coker(i_* F \to i_* F_0)$ can be 
supported only at the points $P_1, \ldots, P_l$, with multiplicities bounded by 
$k_1, \ldots, k_l$, respectively (in the differential geometry picture, cf.
\cite{DK}, $\xi_0$ represents the singular part of a connection which may be 
smoothed out by $F$ but may not acquire any negative coefficients; since the
 multiplicities of $Coker(i_* F \to i_* F_0)$  measure the local change of 
$c_2$ one obtains the bound mentioned). But the proof of Proposition 
\ref{converse} shows that the multulplicities of $Coker(i_* F \to i_* F_0)$  
give an upper bound for the total sum, over all $p$, of similar multiplicites 
for $Coker((i_p)_* F_p \to (i_{p-1})_* F_{p-1})$. Hence the condition
of Corollary \ref{conditions}$(i)$ is rather natural from the point of view
of Uhlenbeck spaces.

\section{Algebraization of principal bundles.} 

Let $G$ be an affine algebraic group over $k$. We keep the notation of Section 
2.1. and consider left principal $G$-bundles which are locally trivial in 
fppf topology. For such a $G$-bundle $P$ (over $\widehat{U}$ or an open
subset $U\subset X$) and any scheme $Y$ over $k$ with left $G$-action, 
denote by $P_Y = G\setminus (Y\times_k P)$ the associated fiber bundle, i.e.
the quotient by the left diagonal action of $G$. For instance, when 
$\rho: G \to H$ is a homomorphism of linear algebraic groups over $k$, we can 
consider a left $G$-action on $H$ given by $g \cdot h = h \rho(g)^{-1}$ and 
then $P_H$ is simply the principal $H$-bundle induced via $\rho$. 

\begin{theorem}
\label{principal} Assume that the
identity component $G^\circ$  is reductive.
Then a principal $G$-bundle $\mathcal{P}$ over the formal scheme
$\widehat{U}$ admits an algebraization if and only if for a
fixed exact representation $G \hookrightarrow GL(V)$ the associated
vector bundle $\mathcal{P}_V$ admits an 
algebraization, i.e. satisfies the conditions of Corollary 
\ref{conditions}.
\end{theorem}
The ``only if" part is obvious. Since by a result of Haboush,
cf. Theorem 3.3 in \cite{Ha}, the quotient $GL(V)/\eta(G)$ is
affine, the ``if" part follows from the following general statement. 
\begin{prop}
\label{reduction}
Let $H$ an affine algebraic group over $k$ and $G$ its closed
subgroup such that $H/G$ is affine. 
Suppose that $\mathcal{P}$ is a principal $G$-bundle
over $\widehat{U}$ such that the associated principal $H$-bundle
$\mathcal{Q} = \mathcal{P}_H$ 
admits an algebraization. Then $\mathcal{P}$ admits an 
algebraization.
\end{prop}
First we  establish a preparatory result. As before, $U\subset X$
is an open subset satisfying $U\cap X_0 = U_0$.
\begin{lemma}
\label{section}
Let $H$ be a linear
algebraic group overe $k$,  $Q$ be a principal $H$-bundle 
on $U$ and $\widehat{Q}$ its completion. 
Let also $Y$ be an affine $H$-variety. Then for any section 
$\widehat{s}: \widehat{U} \to \widehat{Q}_Y$ there exists a
section $s: W\to Q_Y$ on an open subset $W\subset U$ containing 
$U_0$, with completion equal to $\widehat{s}$. If $(W, s)$ and
$(W', s')$ are two such algebraizations, then $s = s'$ on $W\cap W'$.
\end{lemma}
\textit{Proof.} One can find a $H$-invariant linear 
subspace $V^\vee \subset k[Y]$ containing a set of generators
of $k[Y]$ as a $k$-algebra. Then the surjection $Sym^*_k (V^\vee)
\to k[Y]$ gives an $H$-equivariant closed embedding $Y\to V$ into the dual 
space $V$. This induces 
closed embeddings $Q_Y \to Q_V$ and $\widehat{Q}_Y \to \widehat{Q}_V$. 

Therefore $\widehat{s}$ becomes a section of the 
vector bundle $\widehat{Q}_V$. By Proposition \ref{converse}$(ii)$ 
the completion of the coherent sheaf
$i_* Q_V$ is isomorphic to $\widehat{i}_* 
\widehat{Q}_V$ and therefore by Proposition \ref{facts}$(i)$ there exists
a unique section $\tilde{s}$ of $i_* Q_V$ with completion 
given by $\widehat{i}_* \widehat{s}$. Set $s = \tilde{s}|_U$.

It remains to show that $s(W) \subset Q_V$ on some $W$ as above. 
Let $\mathcal{A} = Sym^* (Q_{V}^\vee)$
be the sheaf of symmetric algebras on $U$ corresponding to $Q_V$
and $\mathcal{I}\subset \mathcal{A}$ the ideal sheaf of 
$Q_Y$. The section $s$ gives the evaluation morphism
$\rho: \mathcal{A}\to \mathcal{O}_U$. The sheaf $G = \rho(\mathcal{I})$
is coherent, being a subsheaf of $\mathcal{O}_U$. Since $\widehat{s}$ 
takes values in $\widehat{Q}_Y$, the completion $\widehat{G}$ is zero. 
By Corollary 10.8.12 in \cite{EGAI} this implies 
$Supp(G) \cap U_0 = \emptyset$ hence 
$W = U\setminus Supp(G)$ satisfies the conditions of the lemma. The
uniqueness of $s$ follows from the uniqueness of $\widetilde{s}$.  $\square$

\bigskip
\noindent
\textit{Proof of Proposition \ref{reduction}.} 
Let $(U, Q)$ be an algebraization of $\mathcal{Q}$. 
In general, giving a principal $G$-bundle is equivalent to 
giving a principal $H$-bundle $\mathcal{R}$ 
together with a reduction to $G$, i.e. a section of the 
associated bundle $\mathcal{R}_{H/G}$ with the fiber 
$H/G$.  Since $\mathcal{Q}$ is induced
from $\mathcal{P}$, we get a section $\widehat{s}:\widehat{U}\to 
\mathcal{Q}_{H/G}$ 
and by the above lemma there exists $s: W\to Q_{H/G}$ such that $\widehat{s}$
is equal to its completion. Then $\mathcal{P}$ admits an algebraization
$(W, P)$ where $P$ is the pullback of the principal $G$-bundle 
$Q \to Q_{H/G}$ via $s: W \to Q_{H/G}$. 
$\square$

\section{Categorical formulations.}

\begin{prop} The functor $F \mapsto \widehat{F}|_{\widehat{U}}$ induces
and eqivalence between the full subcategory of all coherent sheaves $E$ on 
$X$ which are locally free at the points of $U_0 \subset X$ and have 
$depth_x E  \geq 2$ at the points where $E$ is not locally free, and the
full subcategory of locally free sheaves on $\widehat{U}$ 
 admitting algebraization.
\end{prop}
\textit{Proof.} Let $(U, F)$ be an algebraization of $\mathcal{F}$.
Then the sheaf $E = i_* F$ satisfies $E \simeq i_* i^* E$ hence
by Proposition \ref{facts}$(iii)$ $depth_x E \geq 2$ for all $x \in Z = 
X \setminus U$. We also observe that $E$ is uniquely determined by 
$\mathcal{F}$,
since by Propositions \ref{facts}$(i)$ and \ref{converse}$(ii)$ it is the 
unique coherent sheaf on $X$ such that $\widehat{E}\simeq \widehat{i}_* 
\mathcal{F}$. Thus the functor described is essentially surjective 
on objects. 
For the morphisms, let $\mathcal{F}_1, 
\mathcal{F}_2$ be a pair of vector bundles on $\widehat{U}$ 
with algebraizations $(U, F_1)$ and $(U, F_2)$, respectively, which
we may assume to be defined on the same $U$. 
Denote by $E_1 = i_*F_1, E_2 = i_*F_2$ the corresponding coherent
sheaves on $X$. Then 
$
Hom_{\widehat{U}}(\mathcal{F}_1, \mathcal{F}_2)
= Hom_{\widehat{X}}(\widehat{i}_* \mathcal{F}_1, \widehat{i}_* \mathcal{F}_2)
= Hom_X (E_1, E_2) 
$
where the first equality is by adjunction of $i^*$ and $i_*$ and the
second by Propositions \ref{facts}$(i)$ and \ref{converse}$(ii)$.  $\square$

\bigskip
\noindent
To formulate a result for principal bundles, let $\mathcal{B}(G, U_0)$
be the groupoid category in which the objects are given by pairs $(U, P)$ where
$U \subset X$ is an open subset with $U\cap X_0 = U_0$, and $P$ is
a principal $G$-bundle on $U$. Morphisms from $(U, P)$ to $(U', P')$ are given
by
the set of equivalence classes of pairs $(W, \psi)$ where $W \subset U \cap U'$
is an open subset with $W\cap X_0 = U_0$ and $\psi: P|_W \to P'|_W$ an
isomorphism of $G$-bundles. Two such pairs $(W, \psi)$ and $(W, \psi')$ are 
equivalent if $\psi = \psi'$ on $W\cap W'$.  Also denote by 
$Bun(G, \widehat{U})$ the groupoid category of $G$-bundles on the formal 
scheme $\widehat{U}$. Completion along $U_0$ defines a functor 
$\Psi: \mathcal{B}(G, U_0) \to Bun(G, \widehat{U})$. 
The following statement summarizes our results on 
algebraization of principal bundles
\begin{theorem} 
\label{categorical}
With the notation of Section 2.1,
\begin{enumerate}
\item For any affine algebraic group $G$ over $k$,
$\Psi: \mathcal{B}(G, U_0)\to Bun(G, \widehat{U})$
is full and strict.
\item For $G = GL_n(k)$ the essential image of $\Psi$ is the full
subcategory of rank $n$ vector bundles $\mathcal{F} = \plimit F_p$ 
on $\widehat{U}$ which satisfy
 the equivalent conditions (i)-(iii) of Corollary 
\ref{conditions}.
\item Let $G \hookrightarrow H$ be a closed embedding of affine algebraic 
groups over $k$
such that $H/G$ is affine. Then the natural functor from $G$-bundles to 
$H$-bundles induces an equivalence of categories
$$
\mathcal{B}(G, U_0) \simeq Bun(G, \widehat{U}) \times_{Bun(H, \widehat{U})}
\mathcal{B}(H, U_0)
$$ 
\end{enumerate}
\end{theorem}
\textit{Proof.} To prove $(i)$ suppose that $\mathcal{P}, \mathcal{P}'$ are two
 principal bundles on $\widehat{U}$ admitting algebraizations $P, P'$, 
respectively, which we may assume to be defined on the same  $U \subset X$.
Let $\widehat{\psi}:
\mathcal{P}\to \mathcal{P'}$ be an isomorphism. We need to prove that 
there exists (perhaps after shrinking $U$) a unique isomorphism  
$\psi: P \to P'$ with completion given by $\widehat{\psi}$.
Let $Isom(P, P')$ be the bundle of isomorphisms $P \to P'$. Considering 
graphs of isomorphisms, we can identify $Isom(P, P') \simeq G \setminus 
(P\times_U P')$. On the other hand, $P\times_U P'$ is a principal 
bundle over $G\times_k G$. Define a left action of $G \times_k G$ on $G$
by $(g, h) \cdot f = g f h^{-1}$, then $G \setminus 
(P\times_U P') \simeq (P \times_U P')_{G}$. 
Since $\widehat{\psi}$ gives a section $\widehat{s}$ of $Isom(\mathcal{P},
\mathcal{P}')$, applying Lemma \ref{section} to $H= G \times_k G$ and 
$Y = G$, we get a unique algebraization $s: W \to (P \times_U P')_G 
\simeq Isom(P, P')|_W$, which corresponds to the required isomorphism $\psi$. 
This proves $(i)$.

The statement of $(ii)$ for objects holds by Corollary \ref{conditions} and 
for morphisms by $(i)$. 

For $(iii)$ first observe that the compositions $\mathcal{B} (G, U_0)
\to \mathcal{B}(H, U_0) \to Bun(H, \widehat{U})$ and $\mathcal{B} (G, U_0)
\to Bun(G, \widehat{U}) \to Bun(H, \widehat{U})$ are canonically isomorphic, 
therefore one does get a functor 
$$
\mathcal{B}(G, U_0) \to Bun(G, \widehat{U}) \times_{Bun(H, \widehat{U})}
\mathcal{B}(H, U_0)
$$
On objects, this functor is an equivalence if for a 
$G$-bundle $\mathcal{P}$ on $\widehat{U}$, an $H$-bundle $Q$ on 
$U \subset X$ and an isomorphism 
$\phi: \mathcal{P}_H \simeq \widehat{Q}$, there
exists an open subset $W \subset U$ with $W\cap X_0 = U_0$, 
a $G$-bundle $P$ on $W$ and isomorphisms $\widehat{P} \simeq
\mathcal{P}$ and $P_H \simeq Q|_W$ which induce $\phi$ in 
a natural way. This is equivalent to finding an algebraization of
the section $\widehat{s}: \widehat{U} \to \widehat{Q}_{H/G}$ induced by $\phi$, 
which was done  in the proof of Proposition \ref{reduction}. 
On  morphisms, without loss of generality it suffices to consider
two $G$-bundles $P, P'$ defined on the same open set $U$, 
and isomorphisms $\psi: P_H \simeq P'_H$, $\widehat{\phi}:
\widehat{P} \to \widehat{P}'$  which have the same image in $Bun(H, 
\widehat{U})$.  We need
to show that there exists a unique isomorphism $\phi: P \to P'$ inducing
$\widehat{\phi}$ and $\psi$ in the natural sense. But by $(i)$ there 
exists a unique $\phi$ with completion equal to $\widehat{\phi}$. 
Since by assumption the isomorphisms $\psi' = \phi_H$ and $\psi$ are
equal after completion, $\psi' = \psi$ by part $(i)$.
This finishes the proof. $\square$

\bigskip
\noindent
\textit{Address:} Department of Mathematics, MSTB 103, UC Irvine, Irvine CA
92697. 

\end{document}